\title{Ambarzumian-type mixed inverse spectral problems for  Jacobi matrices}
\author{Ethan Luo,
        Steven Ning,
        Tarun Rapaka,
        Xuxuan Joyce Zheng
        }
\documentclass{article}
\usepackage{graphicx,amssymb,listings,amsmath}
\usepackage[a4paper, total={6in, 8in}]{geometry}

\begin{document}
\newpage
\maketitle

\begin{abstract}
    We investigate Ambarzumian-type mixed inverse spectral problems for Jacobi matrices. Specifically, we examine whether the Jacobi matrix can be uniquely determined by knowing all but the first $m$ diagonal entries and a set of $m$ ordered  eigenvalues.
\end{abstract}

\section{Introduction}
The Jacobi matrix is a tridiagonal matrix defined as
$$
 \left[
 \begin{array}{rrrrr}
    b_1 & a_1 & 0 & 0 & 0 \\
    a_1 & b_2 & a_2 & \ddots & 0 \\
    0  &  a_2 & b_3 & \ddots & 0 \\
    0 & \ddots & \ddots & \ddots & a_{n-1} \\
    0 & \dots & 0 & a_{n-1} & b_n
 \end{array}
 \right]
$$
where \( n \in \mathbb{N} \), \( b_k \in \mathbb{R} \) for \( k \in \{1, 2, \dots, n\} \), and \( a_k > 0 \) for \( k \in \{1, 2, \dots, n-1\} \). When \( a_k = 1 \) for each \( k \in \{1, 2, \dots, n-1\} \), this matrix defines the finite discrete Schr\"odinger operator, which we refer to as the discrete Schr\"odinger matrix.

In this paper, we study the following problems: Let us denote the \( n \times n \) discrete Schr\"odinger matrix for \( 1 \leq m \leq n \) as follows:

$$
 S_{n,m} := 
 \left[
 \begin{array}{rrrrrr}
    b_1 & 1 & 0 & 0 & \dots & 0 \\
    1 & \ddots & 1 & 0 & \ddots & 0 \\
    0 & 1 & b_m & 1 & \ddots & \vdots \\
    0 & 0 & 1 & 0 & \ddots & 0 \\
    \vdots & \ddots & \ddots & \ddots & \ddots & 1 \\
    0 & \dots & \dots & 0 & 1 & 0
 \end{array}
 \right]
$$
Additionally, let us denote the free discrete Schr\"odinger matrix of size \( n \times n \) by:
$$
 F_{n} := 
 \left[
 \begin{array}{rrrrr}
    0 & 1 & 0 & \dots & 0 \\
    1 & 0 & 1 & \ddots & \vdots \\
    0 & 1 & 0 & \ddots & 0 \\
    \vdots & \ddots & \ddots & \ddots & 1 \\
    0 & \dots & 0 & 1 & 0
 \end{array}
 \right]
$$

In \cite{H2021}, the authors posed the following Ambarzumian-type mixed inverse spectral problem:

\textbf{Q1.} If \( S_{n,m} \) and \( F_n \) share \( m \) consecutive eigenvalues, can we conclude that \( b_1 = \dots = b_m = 0 \) or equivalently \( S_{n,m} = F_n \)?

When $m=1$, Q1 (and Q2) can be affirmatively answered by applying a special case of the result from Gesztesy and Simon \cite[Theorem 4.2 ]{GS1997}.

When \( m = n \), by standard Ambarzumian-type arguments, Q1 can be answered positively \cite [Theorem 3.3]{H2021} (also see \cite{Ambarzumian1929}). When \( m = 2 \), the authors  also answer Q1 positively  \cite[Theorem 4.2]{H2021}. For more history and background related to inverse spectral problems, we refer readers to \cite{H2021}.

In this paper, we ask a more general question:

\textbf{Q2.} If \( S_{n,m} \) and \( F_n \) share \( m \) ordered eigenvalues, can we conclude that \( b_1 = \dots = b_m = 0 \) or that \( S_{n,m} = F_n \)?

Compared to Q1, we do not require the known eigenvalues to be consecutive. Our investigation suggests that for arbitrary \( m \), the answer to Q2 could be negative. Therefore, in this paper, we will focus on the case where \( n = 5 \) and \( m = 3 \). For this scenario, our findings indicate that the answer to Q1 is positive, whereas the answer to Q2 is negative.

\section{Macaulay2}

To get some initial insight on the problem, we wrote the following program (Macaulay2 code can be found at the github repository\footnote{https://eluo5.github.io/inverse-spectral-problem/solver.m2}) for checking all cases of the problem. This resulted in us finding a single case that was not working. For making this program, we first knew that the characteristic polynomial of a matrix $A$ is a polynomial in terms of $\lambda$ defined by $\det(A - \lambda I)$. The equations for this program then came from assigning the characteristic polynomial of $S_{5,3}$ equal to the characteristic polynomial of $F_{5}$. This allowed us to make 5 equations from equating the coefficients of each power of $\lambda$. Each of these were written out by hand and confirmed with the computer.

\begin{lstlisting}
needsPackage "NumericalAlgebraicGeometry"

--- general form
R = CC[a,b,c,d,e,f,g,h]
A = - a - b - c + d + e + f + g + h
B = a*b + a*c + b*c - 4 - d*e - d*f - d*g - d*h - e*f 
    - e*g - e*h - f*g - f*h - g*h
C = 3*a + 2*b + 2*c - a*b*c + d*e*f + d*e*g
    + d*e*h + d*f*g + d*f*h + d*g*h + e*f*g
    + e*f*h + e*g*h + f*g*h
D = -2*a*b - a*c - b*c + 3 - d*e*f*g - d*e*f*h - d*e*g*h
    - d*f*g*h - e*f*g*h
E = a*b*c - a - c + d*e*f*g*h
F = d + sqrt(3) 
G = e + 1
H = f
I = g - 1
J = h - sqrt(3)

--- six cases
K = {A, B, C, D, E, F, G, H} -- works
K = {A, B, C, D, E, F, G, I} -- works
K = {A, B, C, D, E, F, G, J} -- works
K = {A, B, C, D, E, F, H, I} -- fails
K = {A, B, C, D, E, F, H, J} -- works
K = {A, B, C, D, E, G, H, I} -- works

--- run after each case
sol = solveSystem K
\end{lstlisting}

In the code, $a,b,c$ represent the 3 missing matrix values in order, and $d,e,f,g,h$ represent the 5 eigenvalues in order. We then break the problem down into six possible cases, after getting rid of all symmetric cases. Equations $A,B,C,D,E$ are used in every case to equate the coefficients of the two polynomials, and out of equations $F,G,H,I,J$ we pick three to assign specific eigenvalues to be fixed. Thus, we are actually left with 3 unknown matrix values and 2 unknown eigenvalues as expected. In the code we have found that the only questionable case is when we are given the three eigenvalues $-\sqrt{3}$, $0$, and $1$. Once again we have omitted the symmetric cases for simplicity, so we only have one counterexample rather than two. In this case, a counterexample would be the following matrix.

$$
 \left[
 \begin{array}{rrrrr}
    -1.11542462377894 & 1 & 0 & 0 & 0 \\
    1 & 0.527281667822498 & 1 & 0 & 0 \\
    0 & 1 & 0.702345226288011 & 1 & 0 \\
    0 & 0 & 1 & 0 & 1 \\
    0 & 0 & 0 & 1 & 0
 \end{array}
 \right]
$$

The ordered eigenvalues for this matrix are $-\sqrt{3}$, $-1.25874960534751$, $0$, $1$, $2.10500268324796$, approximately. Numerically, this provides a counterexample to  our problem, though we still need to establish a rigorous mathematical proof. The proof appears to be quite challenging, as it involves solving a system of 3 variables and 3 nonlinear equations, which we will demonstrate below. As for all the consecutive eigenvalue cases for Q1, the program verifies that all of these are true.

\section{Specific cases checked by hand}

Suppose that the eigenvalues of $S_{5,3}$ are $\widetilde{\lambda_1} \le \widetilde{\lambda_2}  \le \widetilde{\lambda_3}  \le \widetilde{\lambda_4}  \le \widetilde{\lambda_5}$, and the eigenvalues of $F_5$ are $\lambda_1 = -\sqrt{3}$, $\lambda_2 = -1$, $\lambda_3 = 0$, $\lambda_4 = 1$, and $\lambda_5 = \sqrt{3}$. We then know that the characteristic polynomial of $S_{5,3}$ is $f_{S_{5,3}}(\lambda) = b_1b_2b_3\lambda^2 - b_1b_2b_3 - b_1b_2\lambda^3 + 2b_1b_2\lambda - b_1b_3 \lambda^3 + b_1b_3 \lambda + b_1 \lambda^4 - 3b_1\lambda^2 + b_1 -b_2b_3\lambda^3 + b_2b_3\lambda + b_2\lambda^4-2b_2\lambda^2+b_3\lambda^4-2b_3\lambda^2+b_3-\lambda^5+4\lambda^3-3\lambda$. After some grouping we get that
\begin{equation} \label{eq1}
\begin{split}
f_{S_{5,3}}(\lambda) = &(\lambda^2-1) ((b_1-\lambda)(b_2-\lambda)(b_3-\lambda) - (b_1-\lambda)-(b_3-\lambda)) \\
&+ \lambda((b_1-\lambda)(b_2-\lambda)-1).
\end{split}
\end{equation}

\subsection{Case 1: $\mathbf{\widetilde{\lambda_2} = \lambda_2 = -1}$, $\mathbf{\widetilde{\lambda_3} = \lambda_3 = 0}$, and $\mathbf{\widetilde{\lambda_4} = \lambda_4 = 1}$.}

Since $f_{S_{5,3}}(\widetilde{\lambda_4}) = 0$ and $\widetilde{\lambda_4} = 1$, we get from equation \eqref{eq1} that
\[
1((b_1-1)(b_2-1) - 1) = b_1 b_2 - b_1 - b_2 = 0,
\]
which implies that
$$b_1 b_2 = b_1 + b_2.$$
Similarly since $f_{S_{5,3}}(\widetilde{\lambda_2}) = 0$ and $\widetilde{\lambda_2} = -1$, we get
\[
-1((b_1 + 1)(b_2 + 1) = -(b_1 b_2 + b_1 + b_2) = 0,
\]
which implies that
$$b_1 b_2 = b_1 + b_2 = -(b_1 + b_2).$$
Since $b_1 + b_2 = -(b_1 + b_2),$ we know that $b_1 + b_2 = 0,$ so $b_1 b_2 = 0$ as well. Solving this pair of equations, gives us $b_1 = b_2 = 0.$ Then plugging these values of $b_1$ and $b_2$ back into equation \eqref{eq1}, we get
\[
f_{S_{5,3}}(\lambda) = (\lambda^2-1)(\lambda^2 (b_3-\lambda) + 2\lambda - b_3) + \lambda(\lambda^2-1).
\]
Since $\lambda = 
\widetilde{\lambda_3} = 0$ is also a root of this polynomial, we see that
\[
(0^2-1)(0^2(b_3 - 0) + 2\cdot 0 - b_3) + 0(0^2 - 1) = b_3 = 0.
\]
So $b_1 = b_2 = b_3 = 0$ and $F_5 = S_{5,3}.$

\subsection{Case 2: $\mathbf{\widetilde{\lambda_1} = \lambda_1 = -\sqrt{3}}$, $\mathbf{\widetilde{\lambda_2} = \lambda_2 = -1}$, and  $\mathbf{\widetilde{\lambda_4} = \lambda_4 = 1}$.}

Since $f_{S_{5,3}}(\widetilde{\lambda_4}) = 0$ and $\widetilde{\lambda_4} = 1$, we get from equation \eqref{eq1} that
\[
1((b_1-1)(b_2-1) - 1) = b_1 b_2 - b_1 - b_2 = 0,
\]
which implies that
$$b_1 b_2 = b_1 + b_2.$$
Similarly since $f_{S_{5,3}}(\widetilde{\lambda_2}) = 0$ and $\widetilde{\lambda_2} = -1$, we once again get
\[
-1((b_1 + 1)(b_2 + 1) = -(b_1 b_2 + b_1 + b_2) = 0,
\]
which tells us that
$$b_1 b_2 = b_1 + b_2 = -(b_1 + b_2).$$
Since $b_1 + b_2 = -(b_1 + b_2),$ we know that $b_1 + b_2 = 0,$ so $b_1 b_2 = 0$ as well. Solving this pair of equations, gives us $b_1 = b_2 = 0.$ Then plugging these values of $b_1$ and $b_2$ back into equation \eqref{eq1}, we get
\[
f_{S_{5,3}}(\lambda) = (\lambda^2-1)(\lambda^2 (b_3-\lambda) + 2\lambda - b_3) + \lambda(\lambda^2-1).
\]
Since $\lambda = 
\widetilde{\lambda_1} = -\sqrt{3}$ is also a root of this polynomial, we see that
\[
2(3(b_3+\sqrt{3}) - 2\sqrt{3} - b_3) - \sqrt{3}(2) = 2(2b_3 + \sqrt{3}) - 2\sqrt{3} = 4b_3 = 0,
\]
which once again implies that $b_3 = 0.$ So for this case we also have that $b_1 = b_2 = b_3 = 0$ and $F_5 = S_{5,3}.$

\subsection{Case 3: $\mathbf{\widetilde{\lambda_2} = \lambda_2 = -1}$, $\mathbf{\widetilde{\lambda_4} = \lambda_4 = 1}$, and $\mathbf{\widetilde{\lambda_5} = \lambda_5 = \sqrt{3}}$.}

Similar to Case 1 and Case 2, $\widetilde{\lambda_2} = -1$, and $\widetilde{\lambda_4} = 1$ give us $b_1 b_2 = b_1 + b_2$ and $b_1 b_2 = -(b_1 + b_2)$ respectively. They form a system with the unique solution $b_1 = b_2 = 0$. 

Then since $f_{S_{5,3}}(\widetilde{\lambda_5}) = 0$ and $\widetilde{\lambda_5} = \sqrt{3}$, we get from equation \eqref{eq1}
\begin{equation}\label{eq2}
2b_1b_2b_3 - \sqrt{3}b_1b_2 - 2\sqrt{3}b_1b_3  - 2\sqrt{3}b_2b_3 + b_1 + 3b_2 + 4b_3 = 0.
\end{equation}

Further applying $b_1 = b_2 = 0$ into \eqref{eq2}, we quickly obtain $b_3 = 0$. Therefore, the only possible values of $b_1, b_2, b_3$ such that
$\mathbf{\widetilde{\lambda_2} = \lambda_2  = -1}$, $\mathbf{\widetilde{\lambda_4} = \lambda_4 = 1}$, 
and $\mathbf{\widetilde{\lambda_5} = \lambda_5 = \sqrt{3}}$ are $b_1 = b_2 = b_3 = 0.$ In other words, if these three eigenvalues are equal, then we are forced to have that $F_5 = S_{5,3}.$

\subsection{Case 4: $\mathbf{\widetilde{\lambda_1} = \lambda_1 = -\sqrt{3}}$, $\mathbf{\widetilde{\lambda_3} = \lambda_3 = 0}$, and  $\mathbf{\widetilde{\lambda_4} = \lambda_4 = 1}$}

This is the failed case mentioned above that we will now attempt to enumerate by hand. 
Since $f_{S_{5,3}}(\widetilde{\lambda_4}) = 0$ and $\widetilde{\lambda_4} = 1$, we can follow a similar method to the two cases above to say that $$b_1b_2 = b_1 + b_2.$$ 
Similarly since $f_{S_{5,3}}(\widetilde{\lambda_3}) = 0$ and $\widetilde{\lambda_3} = 0$, we get
$$(-1) ((b_1)(b_2)(b_3) - (b_1) - (b_3)) = 0.$$ 
which tells us that
$$b_1b_2b_3 = b_1+b_3.$$
Then since $f_{S_{5,3}}(\widetilde{\lambda_1}) = 0$ and $\widetilde{\lambda_1} = -\sqrt{3}$, we get after simplification $$2b_1b_2b_3 + \sqrt{3}b_1b_2 + 2\sqrt{3}b_2b_3 + 2\sqrt{3}b_1b_3 + b_1 + 3b_2 + 4b_3 = 0.$$

So now we have 3 variables and 3 equations, and we can once again just plug these into either Macaulay2 or Wolfram Alpha to see that there does exist other solutions. The issue is that the exact forms of these solutions are too complicated to actually work with. For example the exact form of $b_1 \approx -1.11542462377894$ would be $b_1 = \frac{\sqrt{3} - 1}{3 + 3\sqrt{3}} - \frac{10 + 6\sqrt{3}}{(3 + 3\sqrt{3})\sqrt[3]{91 + 51\sqrt{3}-\sqrt{11844 + 6834\sqrt{3}}}} - \frac{\sqrt[3]{91 + 51\sqrt{3}-\sqrt{11844 + 6834\sqrt{3}}}}{3 + 3\sqrt{3}}$. Another way to look at the closed forms of these three variables would be that $b_1$ is one of the six roots of $3b_1^6 - 12b_1^5  + 15b_1^4  - 19b_1^2  + 18b_1  - 6$, $b_2$ is one of the six roots of $b_2^6 - 4b_2^5  + 9b_2^4  - 16b_2^3  + 19b_2^2 - 18b_2 + 6$, and $b_3$ is one of the six roots of $16b_3^6 - 16b_3^5  + 24b_3^4  - 4b_3^3  - 11b_3^2 - 6b_3 + 6$. This solution would act as a counterexample to Q2, saying that it is possible $F_5 \neq S_{5,3}$ even if they share $3$ ordered eigenvalues.

\section{Future work}

In the future, we will provide a rigorous mathematical proof demonstrating that Case 4 in the previous section fails. Specifically, we will show that there exist nontrivial solutions \( b_1, b_2, \) and \( b_3 \) such that the first, second, and fourth ordered eigenvalues of the matrices \( S_{5,3} \) and \( F_5 \) are the same, which answers Q2 negatively.
Additionally, we can investigate why this particular case is the only one that fails and attempt to generalize our findings. We can also scale up to larger matrices to determine if similar patterns persist. Furthermore, there are many other aspects we can examine, such as the number of eigenvalues needed for each missing matrix value and the effects when the missing matrix values are not consecutive. These are all intriguing avenues for further research that may yield significant insights.

\section*{Acknowledgements:}
This project was conducted as part of the  High School Research Program “PReMa” (Program for Research in Mathematics) at Texas A\&M University. 

We would like to thank Dr. Sherry Gong, Dr. Wencai Liu, Dr. Kun Wang (director) and Dr. Zhizhang Xie  for managing the program.
Special thanks to Dr. Wencai Liu for introducing this project and many inspiring discussions.   This work was partially supported by NSF DMS-2015683, DMS-2246031,  DMS-2000345, and the 2024 MAA Dolciani Mathematics Enrichment Grant.

\noindent \textit{Ethan Luo, Texas A\&M University, College Station, Texas 77843. el0337@tamu.edu}

Ethan Luo is currently a senior majoring in computer science and applied mathematics at Texas A\&M. His interests include Graph Theory, Theoretical Computer Science, and Discrete Math. He has worked on several similar research projects using computer algorithms to approximate modern mathematical problems.

Steven Ning is a sophomore from Friendswood High School, with a strong interest in pursuing a variety of fields in mathematics.

Tarun is a high school sophomore from The Woodlands, Texas, who enjoys exploring advanced mathematics.

Xuxuan "Joyce" Zheng is a high school student at College Station High School in Texas, passionate about applied mathematics and its ability to model and solve complex real-world problems.

\end{document}